# Near Time-Optimal Feedback Instantaneous Impact Point (IIP) Guidance Law for Rocket

Byeong-Un Jo[1] and Jaemyung Ahn[2]

*Korea Advanced Institute of Science and Technology (KAIST), 291 Daehak-Ro, Daejeon 34141, Republic of Korea*

## Abstract

This paper proposes a feedback guidance law to move the instantaneous impact point (IIP) of a rocket to a desired location. Analytic expressions relating the time derivatives of an IIP with the external acceleration of the rocket are introduced. A near time-optimal feedback-form guidance law to determine the direction of the acceleration for guiding the IIP is developed using the derivative expressions. The effectiveness of the proposed guidance law, in comparison with the results of open-loop trajectory optimization, was demonstrated through IIP pointing case studies.

## I. Introduction

The instantaneous impact point (IIP) of a rocket, given its position and velocity, is defined as its touchdown point assuming a free-fall flight (without propulsion) [1]. The IIP is considered as a very important information for safe launch operation of a rocket, and it should be calculated and monitored in real-time on the ground facility or on-board of the rocket. The trespassing of IIP trajectory across a destruction line (DL) is one important criterion for a range safety decision – activation of the flight termination system (FTS) – for flight safety operation. A number of studies on prediction of the IIP and their applications for flight safety operation could be found in the literature. These studies include the techniques for computing the IIP in various coordinate systems [2-5], methods on compensation for the effects of gravity perturbation and atmospheric drag [6], expressions for the time derivatives of IIP [7], and introduction of a new flight safety criterion [8].

---
[1] Graduate Research Assistant, Department of Aerospace Engineering, 291 Daehak-Ro.
[2] Associate Professor of Aerospace Engineering, 291 Daehak-Ro; jaemyung.ahn@kaist.ac.kr. (Corresponding Author).

In addition to flight safety operations, the IIP can be used for pre-flight analysis and open-loop optimization of a rocket, particularly to obtain and specify the impact point of separated stages. Yoon and Ahn proposed a trajectory optimization procedure considering the IIPs of the first-stage and payload fairing segments of a launch vehicle as explicit constraints [9]. Using the dispersion analysis, Mandic introduced a guidance and control algorithm that can steer a rocket so that its impact point reaches a target location [10]. The IIP change is important for recent landing guidance of a separated stage of a reusable launch vehicle. For example, it is known that the landing guidance for the separated first stage of Falcon 9 involves the "boostback burn" using three out of nine engines, which change the IIP of stage toward the landing site (barge ship or launch site) [11].

This paper proposes a new near time-optimal feedback guidance law that moves the IIP of a rocket to a target point, whose schematic diagram is shown in Fig. 1. The analytic formulation that describes the time derivatives of an IIP for a given external acceleration vector (primarily produced by the propulsion system) was established. An optimization problem that determines the components of the external acceleration vector to align the IIP derivative vector with the desired direction and maximize its magnitude was formulated, and it can be solved analytically by introducing the Lagrange multipliers. The proposed guidance law was validated through a case study and compared with the results of an open-loop trajectory optimization to minimize the final time.

Three key contributions of this study are summarized as follows. First, the proposed guidance law is a feedback form that can explicitly specify the IIP at the final time. Since the guidance law is a feedback form, it is robust to the error coming from various sources (e.g., the position and velocity errors at the beginning of the guidance). Second, its performance is near time-optimal, and near fuel-optimal assuming the acceleration profile of the rocket is given. Lastly, the proposed law does not involve any iterative procedure, which is a very attractive property for its potential on-board implementation.

This remainder of this paper is organized as follows. Section II introduces the methodology to calculate the Keplerian IIP and its time derivatives. Section III proposes a feedback guidance law to move the IIP of a rocket to a desired point in a near time-optimal manner, which is obtainable by solving a constrained optimization problem established based on the results of Section II. Case studies for validating the guidance law are presented in Section IV. Finally, Section V discusses the comprehensive conclusions of this study and potential opportunities for future work.

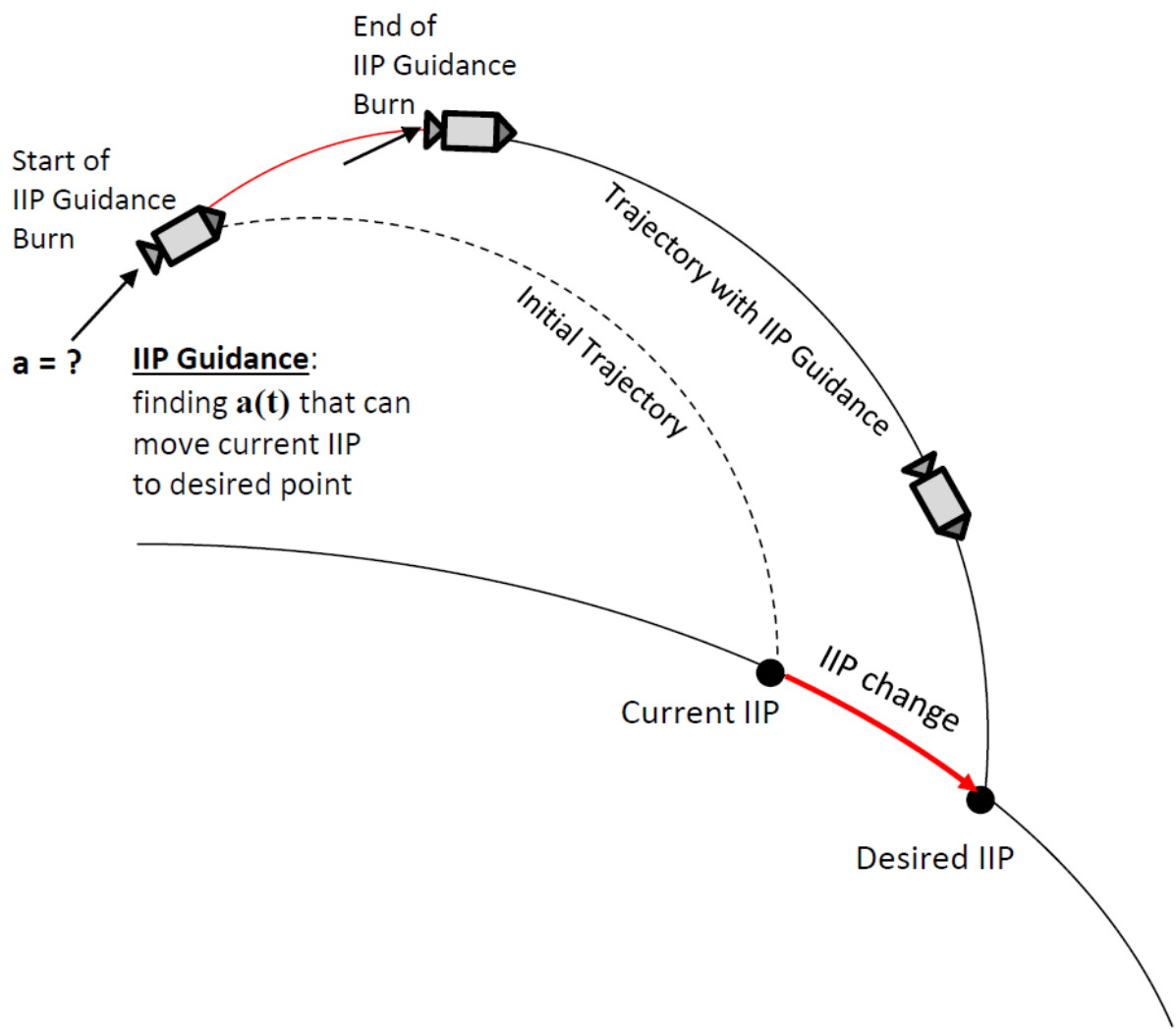

**Figure 1. Schematic diagram explaining IIP guidance**

## II. Calculation of IIP and Its Time Derivatives

This section introduces the procedures to calculate the IIP and its time derivatives in an inertial (Earth centered inertial, ECI) and rotating (Earth centered Earth fixed, ECEF) frames, which provides the fundamentals of the feedback IIP guidance law discussed in this paper. Note that Subsections II-A and II-B are written by summarizing the results of prior studies conducted by Ahn and Roh [5, 7]. The parameters and geometry used to compute the IIP and its time derivatives are shown in Fig. 2.

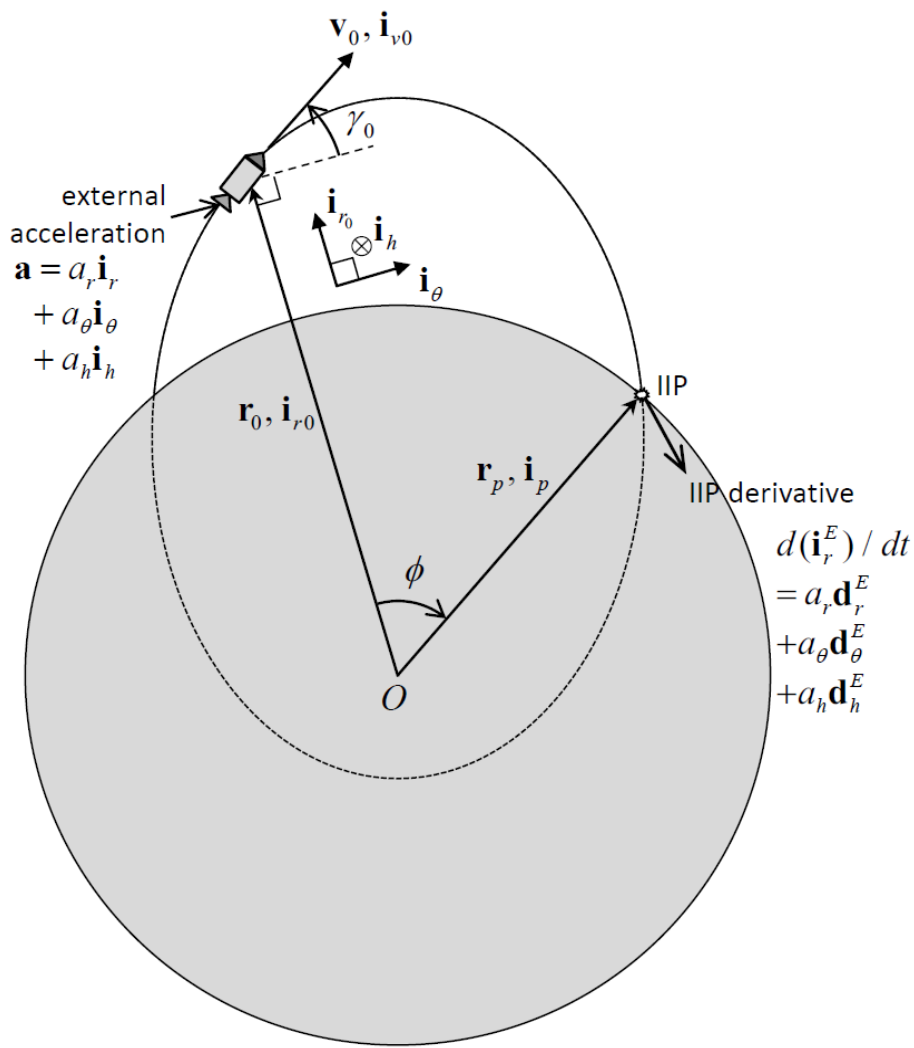

Figure 2. Parameters and geometry for computing IIP and its derivatives

### A. Calculation of Keplerian IIP [5]

Consider the translational motion of a rocket subject to gravity ($\mathbf{g}$) and an external acceleration ($\mathbf{a}$) as follows:

$$\dot{\mathbf{r}} = \mathbf{v} \tag{1}$$

$$\dot{\mathbf{v}} = \mathbf{g} + \mathbf{a} = \mathbf{g}(\mathbf{r}) + a_r \mathbf{i}_r + a_\theta \mathbf{i}_\theta + a_h \mathbf{i}_h \tag{2}$$

In the dynamic equations, $\mathbf{r}$ and $\mathbf{v}$ are position and velocity of the rocket, respectively; $a_r$, $a_\theta$, and $a_h$ are the components of acceleration vector in position, tangential, and linear momentum directions, respectively. In addition, the unit vectors $\mathbf{i}_r$, $\mathbf{i}_\theta$, and $\mathbf{i}_h$ are defined as

$$\mathbf{i}_r \equiv \frac{\mathbf{r}}{\|\mathbf{r}\|} = \frac{\mathbf{r}}{r} \tag{3}$$

$$\mathbf{i}_h \equiv \frac{\mathbf{r} \times \mathbf{v}}{\|\mathbf{r} \times \mathbf{v}\|} = \frac{\mathbf{r} \times \mathbf{v}}{\|\mathbf{h}\|} = \frac{\mathbf{r} \times \mathbf{v}}{h} \tag{4}$$

$$\mathbf{i}_\theta \equiv \mathbf{i}_h \times \mathbf{i}_r \tag{5}$$

If the Keplerian two-body motion is assumed, the gravitational acceleration is expressed as

$$\mathbf{g}(\mathbf{r}) = -\frac{\mu}{\|r\|^3} \mathbf{r} \tag{6}$$

Given current position ($\mathbf{r}_0$) and velocity ($\mathbf{v}_0$) of the rocket, its IIP in the ECI coordinate frame is expressed as follows:

$$\mathbf{i}_p = \frac{\cos(\gamma_0 + \phi)}{\cos \gamma_0} \mathbf{i}_{r0} + \frac{\sin \phi}{\cos \gamma_0} \mathbf{i}_{v0} \tag{7}$$

In this equation, $\gamma_0$ and $\phi$ are respectively the flight path angle and the angle of flight of the rocket expressed as

$$\gamma_0 = \sin^{-1}\left(\frac{\mathbf{r}_0 \cdot \mathbf{v}_0}{\|\mathbf{r}_0\|\|\mathbf{v}_0\|}\right) = \sin^{-1}\left(\frac{\mathbf{r}_0 \cdot \mathbf{v}_0}{r_0 v_0}\right) \tag{8}$$

$$\phi = \sin^{-1}\left(\frac{c_1 c_3 + \sqrt{c_1^2 c_3^2 - (c_1^2 + c_2^2)(c_3^2 - c_2^2)}}{c_1^2 + c_2^2}\right) \tag{9}$$

where $c_1$, $c_2$, and $c_3$ are expressed as

$$c_1 = -\frac{h}{\mu r_0}(\mathbf{r}_0 \cdot \mathbf{v}_0), \ c_2 = \frac{h^2}{\mu r_0} - 1, \ c_3 = \frac{h^2}{\mu r_p} - 1 \tag{10}$$

The time of flight of the launch vehicle – between the current time and the impact time – is expressed as follows [12]

$$t_F = \frac{r_0}{v_0 \cos\gamma_0} \left( \frac{\tan\gamma_0(1-\cos\phi) + (1-\Lambda)\sin\phi}{(2-\Lambda)\left(\frac{1-\cos\phi}{\Lambda\cos^2\gamma_0} + \frac{\cos(\gamma_0+\phi)}{\cos\gamma_0}\right)} + \frac{2\cos\gamma_0}{\Lambda\left(\frac{2}{\Lambda}-1\right)} \tan^{-1}\left( \frac{\sqrt{\frac{2}{\Lambda}-1}}{\cos\gamma_0 \cot\left(\frac{\phi}{2}\right) - \sin\gamma_0} \right) \right) \quad (11)$$

where $\Lambda \, (\equiv (v_0/v_c)^2 = r_0 v_0^2/\mu)$ is defined as the square of the ratio between the current velocity and the circular orbit velocity with given radius ($v_c = \sqrt{\mu/r_0}$). The IIP latitude and longitude in the ECI coordinate system can be expressed using the components of the IIP unit vector in Eq. (7) as

$$\text{Lat}_p = \sin^{-1}(i_{pz}) \quad (12)$$

$$\text{Lon}_p = \tan^{-1}(i_{py}, i_{px}) \quad (13)$$

The IIP longitude in the ECEF coordinate system is obtained by reflecting the Earth's rotation during the time of flight as

$$\text{Lon}_p^E = \text{Lon}_p - \omega_e(t - t_{ref} + t_F) = \text{Lon}_p - \omega_e \Delta t \quad (14)$$

where $\omega_e$ indicates the rotational rate of the Earth, $t$ is the current time, and $t_{ref}$ is the time when the ECI and ECEF coordinates coincide. For details on the procedure to compute the flight time, refer to references [5, 13].

## B. Time Derivatives of Keplerian IIP [7]

Change in the IIP of a rocket occurs when the external acceleration is not zero ($\mathbf{a} \neq \mathbf{0}$). The time derivative of the IIP unit vector ($d(\mathbf{i}_p)/dt$) is expressed as a linear combination of components of the external acceleration vector as follows:

$$\frac{d}{dt}\mathbf{i}_p = \begin{bmatrix} \frac{di_{px}}{dt} & \frac{di_{py}}{dt} & \frac{di_{pz}}{dt} \end{bmatrix}^T = a_r \mathbf{d}_r + a_\theta \mathbf{d}_\theta + a_h \mathbf{d}_h \quad (15)$$

In this equation, $\mathbf{d}_r, \mathbf{d}_\theta$, and $\mathbf{d}_h$ are the directions of the IIP derivative associated with three components of external acceleration ($a_r, a_\theta$, and $a_h$), defined as follows:

$$\mathbf{d}_r = \frac{D_{a_r}^{\dot\phi}}{h}\left[-(h\sin\phi + (\mathbf{r}_0\cdot\mathbf{v}_0)\cos\phi)\cdot\mathbf{i}_{r0} + (r_0 v_0 \cos\phi)\cdot\mathbf{i}_{v0}\right] \quad (16)$$

$$\mathbf{d}_\theta = \frac{D_{a_\theta}^{\dot\phi}}{h}\left[-\left(h\sin\phi + (\mathbf{r}_0\cdot\mathbf{v}_0)\cos\phi\right)\cdot\mathbf{i}_{r0} + (r_0 v_0 \cos\phi)\cdot\mathbf{i}_{v0}\right] \quad (17)$$

$$\mathbf{d}_h = \frac{1}{h}(r_0 \sin\phi)\cdot\mathbf{i}_h \quad (18)$$

In Eqs. (16)-(18), $D_{a_r}^{\dot\phi}$ and $D_{a_\theta}^{\dot\phi}$ represent the influence of $a_r$ and $a_\theta$ on $\dot\phi$, and are expressed as

$$D_{a_r}^{\dot\phi} = \frac{\partial\dot\phi}{\partial a_r} = \frac{h\sin\phi}{\mu(-c_2\sin\phi + c_1\cos\phi)} \quad (19)$$

$$D_{a_\theta}^{\dot\phi} = \frac{\partial\dot\phi}{\partial a_\theta} = \frac{2h(r_0/r_e - \cos\phi) + (\mathbf{r}_0\cdot\mathbf{v}_0)\sin\phi}{\mu(-c_2\sin\phi + c_1\cos\phi)} \quad (20)$$

Note that, owing to the property that IIP moves along the surface of the Earth, the time derivative components $\mathbf{d}_r$, $\mathbf{d}_\theta$, and $\mathbf{d}_h$ are all tangential to $\mathbf{i}_p$. In addition, Eqs. (16) and (17) indicate that $\mathbf{d}_r$ and $\mathbf{d}_\theta$, which are the in-plane acceleration components associated with $a_r$ and $a_\theta$, are parallel. On the other hand, $a_h$ generates the plane change motion of the rocket, whose contribution in the IIP derivative is directed to $\mathbf{d}_h$.

The time derivative of the time of flight ($\dot t_F$) is expressed as follows

$$\dot t_F = -1 + D_{a_r}^{t_F} a_r + D_{a_\theta}^{t_F} a_\theta = -1 + \left(D_{\dot a}^{t_F} D_{a_r}^{\dot a} + D_{\dot e}^{t_F} D_{a_r}^{\dot e}\right)a_r^T + \left(D_{\dot a}^{t_F} D_{a_\theta}^{\dot a} + D_{\dot e}^{t_F} D_{a_\theta}^{\dot e}\right)a_\theta \equiv -1 + \dot t_{F,a} \quad (21)$$

where $D_{\dot a}^{t_F}$, $D_{\dot e}^{t_F}$, $D_{a_r}^{\dot a}$, $D_{a_r}^{\dot e}$, $D_{a_\theta}^{\dot a}$, and $D_{a_\theta}^{\dot e}$ in Eq. (21) are given as

$$D_{\dot a}^{t_F} = \frac{\partial t_F}{\partial \dot a} = \frac{3t_F}{2a} - \frac{1}{a^2 en}\left[\frac{r_e(1-e\cos E_p)}{\sin E_p} - \frac{r_0(1-e\cos E_0)}{\sin E_0}\right] \quad (22)$$

$$D_{\dot{e}}^{\dot{t}_F} = \frac{\partial \dot{t}_F}{\partial \dot{e}} = \frac{1}{n}\left[\left(\frac{\cos E_p (1-e\cos E_p)}{e\sin E_p} - \frac{\cos E_0 (1-e\cos E_0)}{e\sin E_0}\right) - (\sin E_p - \sin E_0)\right] \quad (23)$$

$$D_{a_r}^{\dot{a}} = \frac{2a^2 v_0 \sin\gamma_0}{\mu}, \quad D_{a_r}^{\dot{e}} = \frac{pv_0 \sin\gamma_0}{\mu e} \quad (24)$$

$$D_{a_\theta}^{\dot{a}} = \frac{2a^2 v_0 \cos\gamma_0}{\mu}, \quad D_{a_\theta}^{\dot{e}} = \frac{(pa - r_0^2) v_0 \cos\gamma_0}{\mu a e} \quad (25)$$

In Eqs. (24) and (25), $a$, $e$, $p$, and $n$ are the semimajor axis, eccentricity, semiparameter, and mean motion of the orbit, respectively, and $E_0$ and $E_p$ are the eccentric anomaly values of the rocket at the current and impact points expressed as follows

$$a = \frac{1}{1/r_0 - v_0^2/\mu} \quad (26)$$

$$p = \frac{h^2}{\mu} \quad (27)$$

$$n = \sqrt{\mu/a^3} \quad (28)$$

$$e = \sqrt{1 - p/a} \quad (29)$$

$$E_0 = \cos^{-1}\left(\frac{a - r_0}{ae}\right) \quad (30)$$

$$E_p = \cos^{-1}\left(\frac{a - r_p}{ae}\right) \quad (31)$$

Note that $\dot{t}_F$ is composed of two parts: the part caused by gravity (-1) and the derivative created by external acceleration ($\dot{t}_{F,a}$). The second part ($\dot{t}_{F,a}$) is used to define the time derivative of IIP in a rotating (ECEF) frame, which is discussed in the next subsection.

## C. Time Derivatives of Keplerian IIP in ECEF Frame

The time derivative of the Keplerian IIP in the ECEF coordinate frame is expressed as follows

$$\frac{d}{dt}\mathbf{i}_p^R = \mathbf{T}_I^E (\frac{d}{dt}\mathbf{i}_p) - \dot{t}_{F,a}(\Omega_e \times \mathbf{i}_p^E) \tag{32}$$

where $\mathbf{i}_p^E$ is the IIP unit vector expressed in the ECEF frame, $\mathbf{T}_I^E$ is the transform matrix from the ECI to ECEF coordinate system, and $\Omega_e$ is the vector representing the rotation of the Earth expressed as

$$\mathbf{i}_p^E = [i_{px}^E \; i_{py}^E \; i_{pz}^E]^T = [\cos(\text{Lat}_p)\cos(\text{Lon}_p^E), \; \cos(\text{Lat}_p)\sin(\text{Lon}_p^E), \; \sin(\text{Lat}_p)]^T \tag{33}$$

$$\mathbf{T}_I^E = \begin{bmatrix} \cos(\omega_e \Delta t) & \sin(\omega_e \Delta t) & 0 \\ -\sin(\omega_e \Delta t) & \cos(\omega_e \Delta t) & 0 \\ 0 & 0 & 1 \end{bmatrix} \tag{34}$$

$$\Omega_e = [0, \; 0, \; \omega_e]^T \tag{35}$$

Combining Eqs. (15), (21), and (32), one can obtain the following expression for time derivative of IIP in ECEF frame.

$$\frac{d}{dt}\mathbf{i}_p^E = \mathbf{d}_r^E a_r + \mathbf{d}_\theta^E a_\theta + \mathbf{d}_h^E a_h \tag{36}$$

where

$$\mathbf{d}_r^E = \left( \omega_e D_{a_r}^{i_F} \begin{bmatrix} i_{py}^E \\ -i_{px}^E \\ 0 \end{bmatrix} + \mathbf{T}_I^E \mathbf{d}_r \right), \tag{37}$$

$$\mathbf{d}_\theta^E = \left( \omega_e D_{a_\theta}^{i_F} \begin{bmatrix} i_{py}^E \\ -i_{px}^E \\ 0 \end{bmatrix} + \mathbf{T}_I^E \mathbf{d}_\theta \right), \tag{38}$$

$$\mathbf{d}_h^E = \mathbf{T}_I^E \mathbf{d}_h \tag{39}$$

## III. Feedback Law for Near Time-Optimal IIP Guidance

This section proposes a near time-optimal feedback guidance law for changing the current IIP to the desired location. The guidance law ensures that the IIP rate vector is aligned with the direction of a circular arc departing from the current IIP and arriving at the target position, with its maximum magnitude possible.

As represented in Fig. 3, the shortest path between the current and desired IIP along the surface of the Earth is the arc connecting the two points. Let $\mathbf{i}_p^{t,E}$ denote the target IIP vector in ECEF coordinate system. Then, $\mathbf{q}^E$ is defined as a vector normal to the plane specified by $\mathbf{i}_p^{t,E}$ and $(\mathbf{i}_p^{t,E} - \mathbf{i}_p^E)$, and $\mathbf{i}_q^E$ is defined as its unit vector as follows:

$$\mathbf{q}^E \equiv \mathbf{i}_p^E \times (\mathbf{i}_p^{t,E} - \mathbf{i}_p^E) \tag{40}$$

$$\mathbf{i}_q^E = \mathbf{q}^E / \|\mathbf{q}^E\| \tag{41}$$

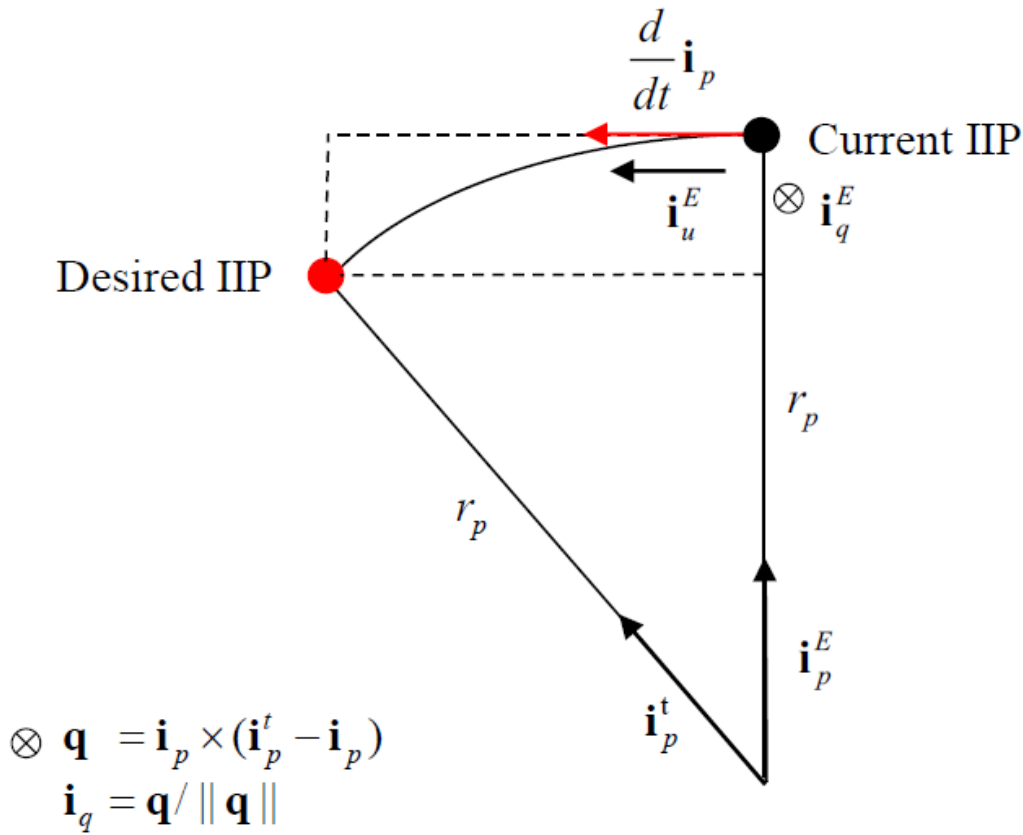

Figure 3: Definitions of unit vectors used for IIP guidance command generation

The unit vector in the direction of the arc (shortest path) between the current and the target IIPs ($\mathbf{i}_u^E$) is defined as follows.

$$\mathbf{i}_u^E = \mathbf{i}_q^E \times \mathbf{i}_p^E \tag{42}$$

The problem determines the acceleration vector $\mathbf{a} = [a_r \ a_\theta \ a_h]^T$ whose magnitude is given as $a_m$. The IIP derivative owing to $\mathbf{a}$ ($= d(\mathbf{i}_p^E)/dt = a_r \mathbf{d}_r^E + a_\theta \mathbf{d}_\theta^E + a_h \mathbf{d}_h^E$) should be parallel to $\mathbf{i}_u^E$ and its magnitude should be maximized. That is, we want to maximize $[d(\mathbf{i}_p^E)/dt] \cdot \mathbf{i}_u^E$ subject to two constraints: 1) $[d(\mathbf{i}_p^E)/dt] \cdot \mathbf{i}_q^E = 0$ and 2) $\mathbf{a} \cdot \mathbf{a} - a_m^2 = 0$. This problem (**P**) can be formulated as a constrained optimization problem described as follows.

[**P**$_{CG}$: IIP guidance acceleration generation]

$$\max_{\mathbf{x}} J(\mathbf{x}) = (\mathbf{c}^T \mathbf{x}) = c_1 x_1 + c_2 x_2 + c_3 x_3 \tag{43}$$

subject to

$$\mathbf{x}^T \mathbf{x} - a_m^2 = x_1^2 + x_2^2 + x_3^2 - a_m^2 = 0 \tag{44}$$

$$\mathbf{f}^T \mathbf{x} = f_1 x_1 + f_2 x_2 + f_3 x_3 = 0 \tag{45}$$

where $\mathbf{x}$ is the decision vector and $\mathbf{c}^T$ and $\mathbf{f}^T$ are the problem parameters defining the objective function and constraint, respectively, defined as

$$\mathbf{x} = [a_r \ a_\theta \ a_h]^T \tag{46}$$

$$\mathbf{c}^T = [\mathbf{i}_u^E \cdot \mathbf{d}_r^E, \ \mathbf{i}_u^E \cdot \mathbf{d}_\theta^E, \ \mathbf{i}_u^E \cdot \mathbf{d}_h^E] = [c_1, c_2, c_3] \tag{47}$$

$$\mathbf{f}^T = [\mathbf{i}_q^E \cdot \mathbf{d}_r^E, \ \mathbf{i}_q^E \cdot \mathbf{d}_\theta^E, \ \mathbf{i}_q^E \cdot \mathbf{d}_h^E] = [f_1, f_2, f_3] \tag{48}$$

The constrained optimization problem defined by Eqs. (43)-(48) can be solved by introducing Lagrange multipliers associated with Eqs (44) and (45), $\boldsymbol{\lambda} (\equiv [\lambda_1, \lambda_2]^T)$. The augmented cost function $J_a(\mathbf{x}, \boldsymbol{\lambda})$ can be defined as

$$\begin{aligned} J_a(\mathbf{x}, \boldsymbol{\lambda}) &= \mathbf{c}^T \mathbf{x} + \lambda_1 (\mathbf{x}^T \mathbf{x} - a_m^2) + \lambda_2 (\mathbf{f}^T \mathbf{x}) \\ &= (c_1 x_1 + c_2 x_2 + c_3 x_3) + \lambda_1 (x_1^2 + x_2^2 + x_3^2 - a_m^2) + \lambda_2 (f_1 x_1 + f_2 x_2 + f_3 x_3) \end{aligned} \tag{49}$$

By applying $\partial J_a / \partial x_i = 0$ (for i = 1, 2, 3), the following equations are obtained.

$$\frac{\partial J_a}{\partial x_1} = c_1 + 2x_1\lambda_1 + f_1\lambda_2 = 0, \quad x_1 = \frac{(-f_1\lambda_2 - c_1)}{2\lambda_1} \tag{50}$$

$$\frac{\partial J_a}{\partial x_2} = c_2 + 2x_2\lambda_1 + f_2\lambda_2 = 0, \quad x_2 = \frac{(-f_2\lambda_2 - c_2)}{2\lambda_1} \tag{51}$$

$$\frac{\partial J_a}{\partial x_3} = c_3 + 2x_3\lambda_1 + f_3\lambda_2 = 0, \quad x_3 = \frac{(-f_3\lambda_2 - c_3)}{2\lambda_1} \tag{52}$$

Putting Eqs. (50)-(52) into Eq. (45) and solving for the Lagrange multiplier $\lambda_2$ yields

$$\lambda_2 = \frac{-\mathbf{c}^T\mathbf{f}}{\mathbf{f}^T\mathbf{f}} = \frac{-(c_1f_1 + c_2f_2 + c_3f_3)}{f_1^2 + f_2^2 + f_3^2} \tag{53}$$

In addition, by combining Eqs. (50)-(52) with Eq. (44), $\lambda_1$ is determined as

$$\lambda_1 = \begin{cases} \dfrac{\sqrt{(-f_1\lambda_2 - c_1)^2 + (-f_2\lambda_2 - c_2)^2 + (-f_3\lambda_2 - c_3)^2}}{2a_m} & (if\ -\mathbf{c}^T(\lambda_2\mathbf{f} + \mathbf{c}) > 0) \\ \dfrac{-\sqrt{(-f_1\lambda_2 - c_1)^2 + (-f_2\lambda_2 - c_2)^2 + (-f_3\lambda_2 - c_3)^2}}{2a_m} & (if\ -\mathbf{c}^T(\lambda_2\mathbf{f} + \mathbf{c}) \leq 0) \end{cases} \tag{54}$$

The acceleration components are given as

$$a_r = x_1 = \frac{-f_1\lambda_2 - c_1}{2\lambda_1} \tag{55}$$

$$a_\theta = x_2 = \frac{-f_2\lambda_2 - c_2}{2\lambda_1} \tag{56}$$

$$a_h = x_3 = \frac{-f_3\lambda_2 - c_3}{2\lambda_1} \tag{57}$$

The inputs to the proposed guidance law are the current position ($\mathbf{r}_0$) and velocity ($\mathbf{v}_0$), which are used to calculate the current IIP and directions for IIP derivatives, and the outputs of the guidance law are the components of the guidance acceleration vector ($a_r, a_\theta, a_h$). In addition, it should be noted that the procedure does not involve any iteration loops, which is an important property for in-flight implementation.

The overall structure of the feedback IIP guidance law presented in this section is shown as a block diagram (with related equation numbers) in Fig. 4. The proposed guidance is a feedback form that uses the position (**r**) and velocity (**v**) as the inputs, computes the IIP ($\mathbf{i}_p^E$) and its direction vectors ($\mathbf{d}_r^E, \mathbf{d}_\theta^E, \mathbf{d}_h^E$) as intermediate parameters, and generates the acceleration command (**a**) as the output. In addition, all the procedures are sequential without any involvement of iterations; even the optimal command generation problem ($\mathbf{P}_{CG}$) is solved analytically.

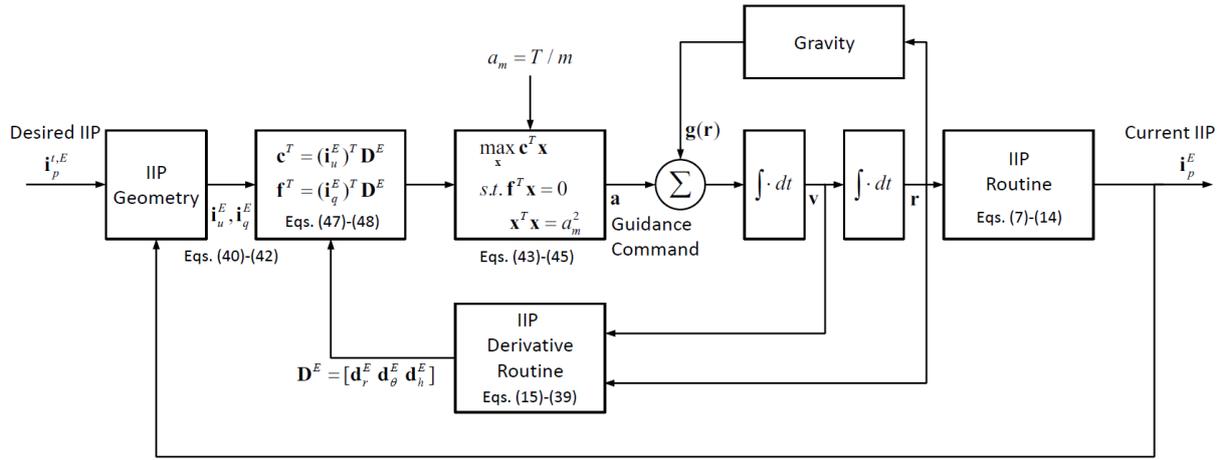

**Figure 4. Block diagram describing the proposed feedback IIP guidance procedure**

## IV. Case Study

The validity of the proposed IIP guidance algorithm presented in Section III is demonstrated through a case study. The IIP change maneuver of the separated first stage of a reusable launch vehicle (Falcon 9 of Space X) was selected as the scenario for the case study. Table 1 summarizes the configuration of the rocket (separated stage) and the initial condition used for the case study.

**Table 1. Rocket configuration and initial condition for case study**

| Rocket Configuration Parameter | Value |
| --- | --- |
| Dry mass, ton | 22.2 |
| Propellant mass, ton | 57 |
| Thrust, tonf | 279.6 |
| Specific impulse, s | 311 |
| **Flight Condition** | **Value** |
| Initial altitude, km | 125.5 |
| Initial speed – inertial ($v_0$), km/s | 1.843 |
| Initial flightpath angle ($\gamma_0$), deg | 3.9 |
| Initial position ($\mathbf{r}_0$), [km, km, km] | [1164, -5507, 3258] |
| Initial velocity ($\mathbf{v}_0$), [m/s, m/s, m/s] | [1337, 743, 1029] |
| Reference time ($t_{ref}$), s | -240 |

Fig. 5 shows the current location, the original IIP, and the desired impact location for the case study. The first and second cases decrease the IIP distance along the downrange direction by 100 km / 200 km. The third and fourth cases increase the ranges by 100 km / 200 km. The fifth case changes the IIP to the crossrange direction by 150 km.

To demonstrate the quality of the solution, the results of the feedback IIP guidance simulation were compared with the open-loop trajectory optimization. The *General Pseudospectral Optimal Control Software* (*GPOPS*) developed by Rao et al. [14] was used to create the open-loop optimal solutions for cases. The objective of the open-loop trajectory optimization was chosen as minimization of the final time ($t_f$), the direction of the acceleration were used and the control variables, and the IIP at $t_f$ was imposed as a constraint ($\mathbf{i}_p^{t,E} = \mathbf{i}_p^E(t_f)$). Both the IIP guidance simulation and open-loop optimization were implemented using MATLAB R2015b and run on a machine with a 3.5 GHz Intel Core i7 CPU, 16 GB RAM, and Windows 7 operating system.

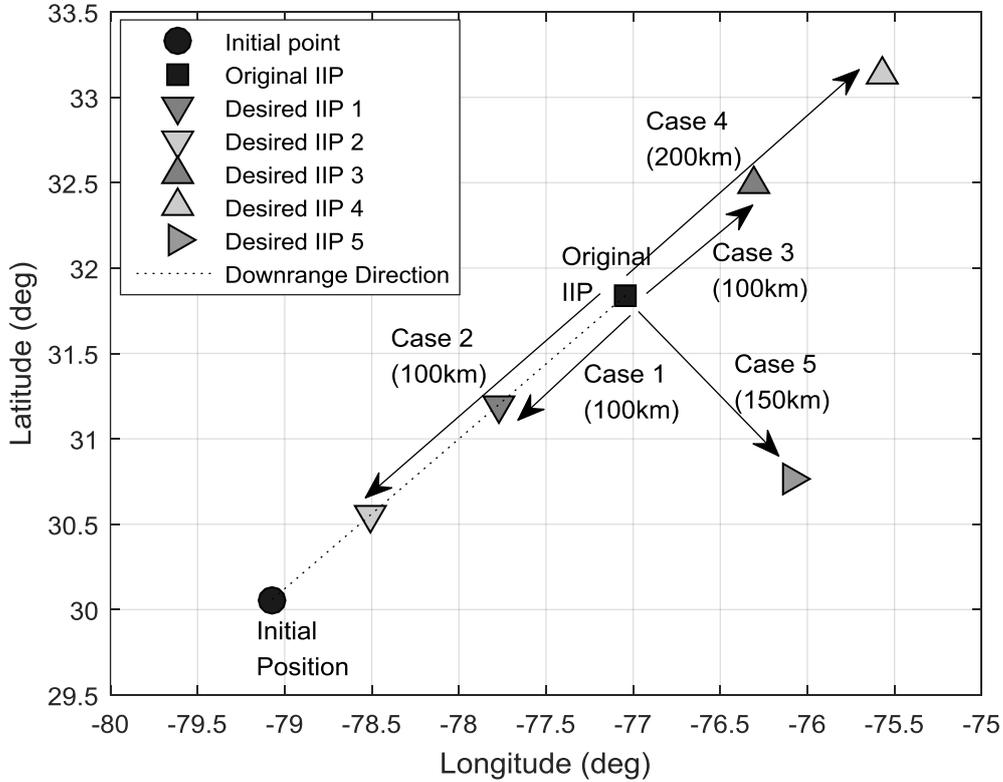

**Figure 5: Current and impact locations for case study**

Table 2 summarizes the results of the case study obtained by using the proposed IIP guidance law and the open-loop trajectory optimization. It can be seen that the differences in the objective function (final time) between the results obtained by two methods are very small – less than 1 % for three cases, and less than 5 % even for the worst case (Case 2). When the results are interpreted as the values of the velocity increments during the maneuver, the performance difference was just 6.1 % at the worst case. Based on this comparison results, we can conclude that the proposed IIP guidance law is near time-optimal with relatively small amount of optimality gap. It is observed that the performance gap (between the IIP guidance and the optimal results) is relatively large when the change in IIP is not aligned with the downrange direction. Case 2 requires that the IIP change to the opposite of the downrange direction and Case 5 changes the IIP in a direction perpendicular to the downrange direction.

**Table 2. Case study results (comparison with open-loop optimization)**

| Case | Final Time ($t_f$), s | | | Propellant Consumption, ton | | $\Delta V$, m/s | | |
|---|---|---|---|---|---|---|---|---|
| | Optimization | IIP Guidance | Difference, % | Optimization | IIP Guidance | Optimization | IIP Guidance | Difference, % |
| 1 | 12.4 | 12.5 | 0.8 | 11.2 | 11.2 | 464 | 465 | 0.3 |
| 2 | 28.1 | 29.5 | 5.0 | 25.2 | 26.5 | 1169 | 1247 | 6.1 |
| 3 | 9.4 | 9.4 | 0.0 | 8.4 | 8.4 | 342 | 343 | 0.0 |
| 4 | 16.4 | 16.4 | 0.0 | 14.7 | 14.7 | 627 | 627 | 0.1 |
| 5 | 21.3 | 21.8 | 2.3 | 19.2 | 19.6 | 846 | 868 | 2.7 |

Figs. 6-8 show the ground trajectories of the rocket position and IIP obtained by the two methods for Cases 2, 4, and 5, where the IIPs are moving toward the desired locations successfully. There is little noticeable difference between the two trajectories.

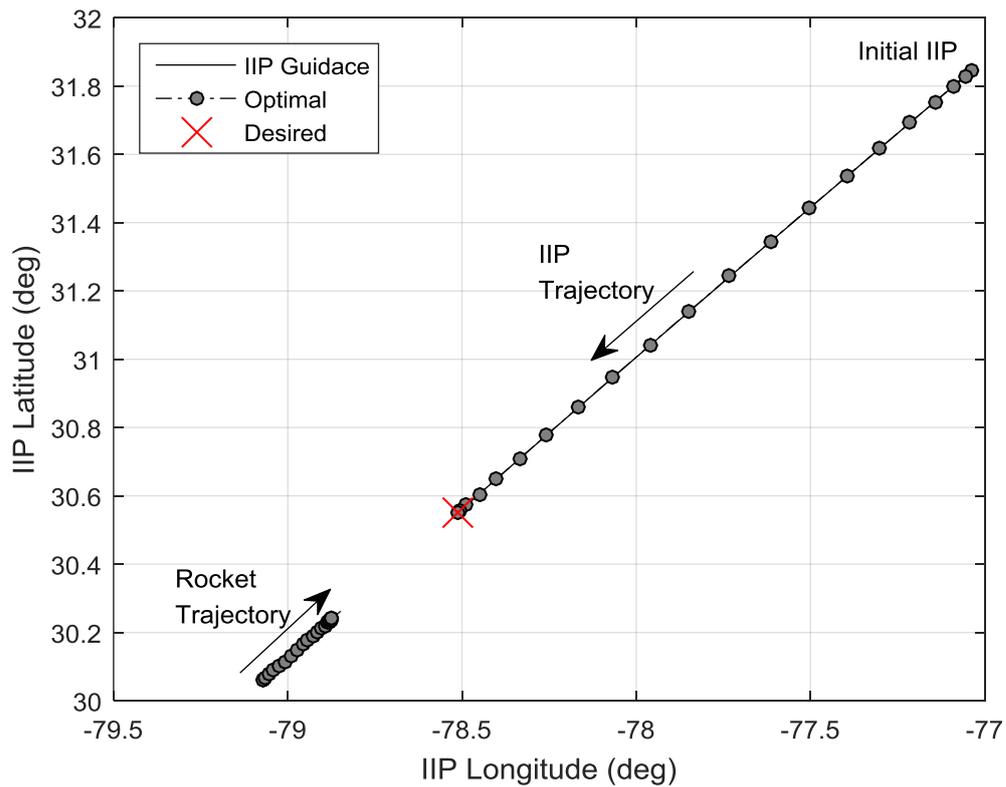

**Figure 6. Trajectories of rocket and its IIP obtained by two methods (Case 2)**

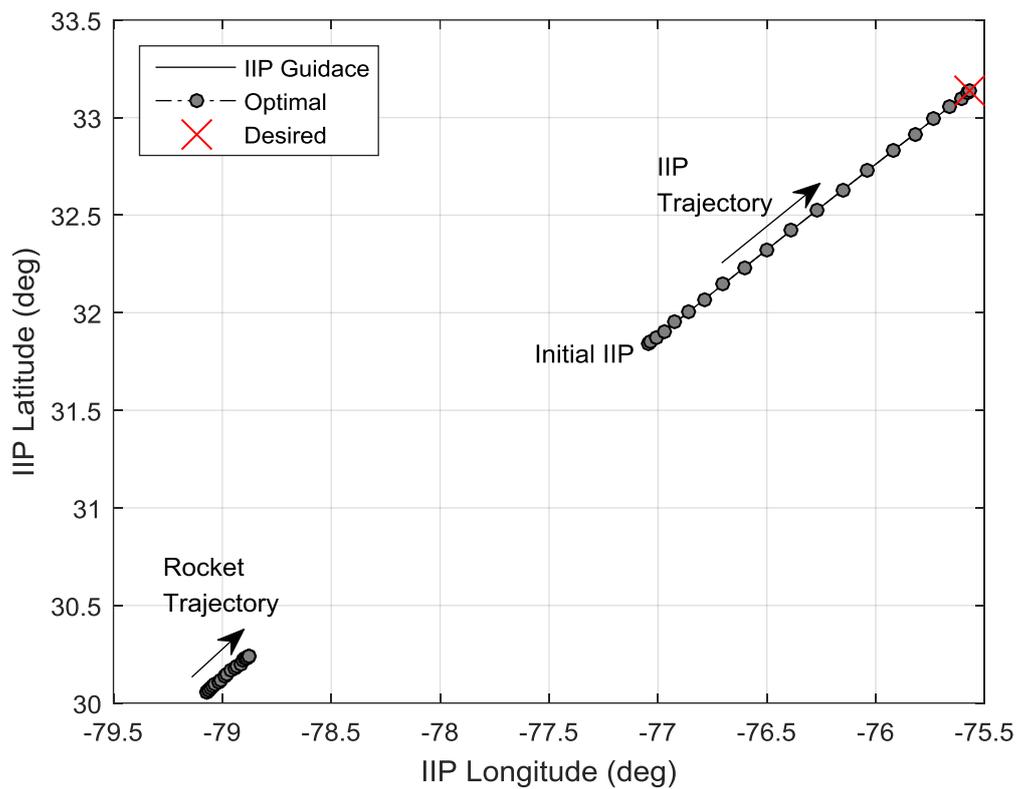

**Figure 7. Trajectories of rocket and its IIP obtained by two methods (Case 4)**

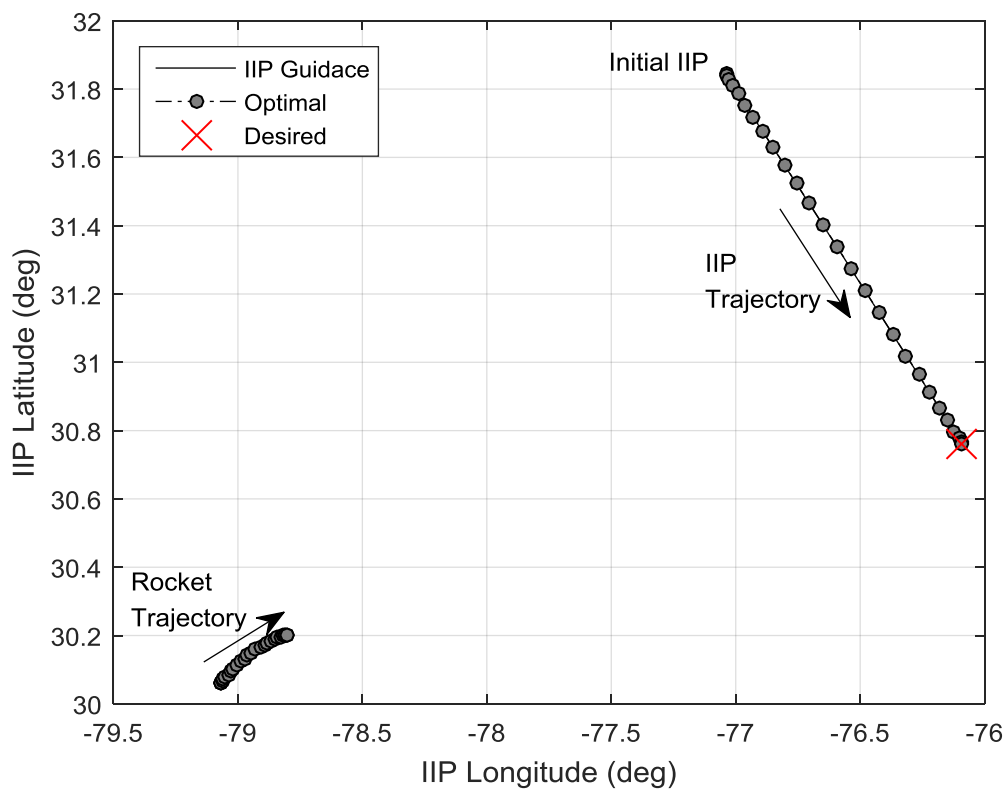

**Figure 8. Trajectories of rocket and its IIP obtained by two methods (Case 5)**

The trajectories using the two methods almost perfectly match. However, the acceleration profiles have meaningful differences, particularly in Cases 2 and 5. Figs. 9-11 compare the histories of acceleration components ($a_r$, $a_\theta$, $a_h$) obtained by the IIP guidance and the open-loop optimization for Cases 2, 4, and 5. The differences in acceleration components for Case 4 presented in Fig. 10, which moves the IIP in the downrange direction, were relatively low. On the contrary, the profiles of acceleration components for the cases where the direction of IIP change is significantly different from the downrange direction (Case 2: the opposite direction, Case 5: 90 degree difference). It was observed that the differences in components associated with in-plane maneuver ($a_r, a_\theta$) are relatively large compared with the difference in $a_h$, which governs the plane change maneuver.

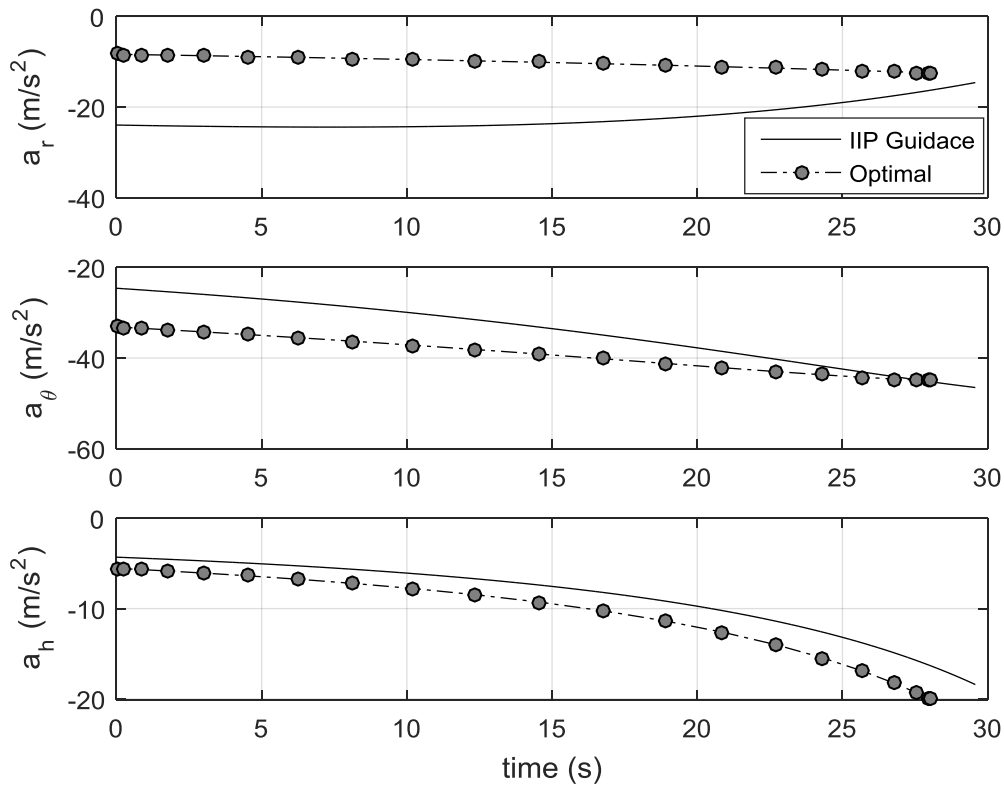

**Figure 9: Histories of acceleration components obtained by two methods (Case 2)**

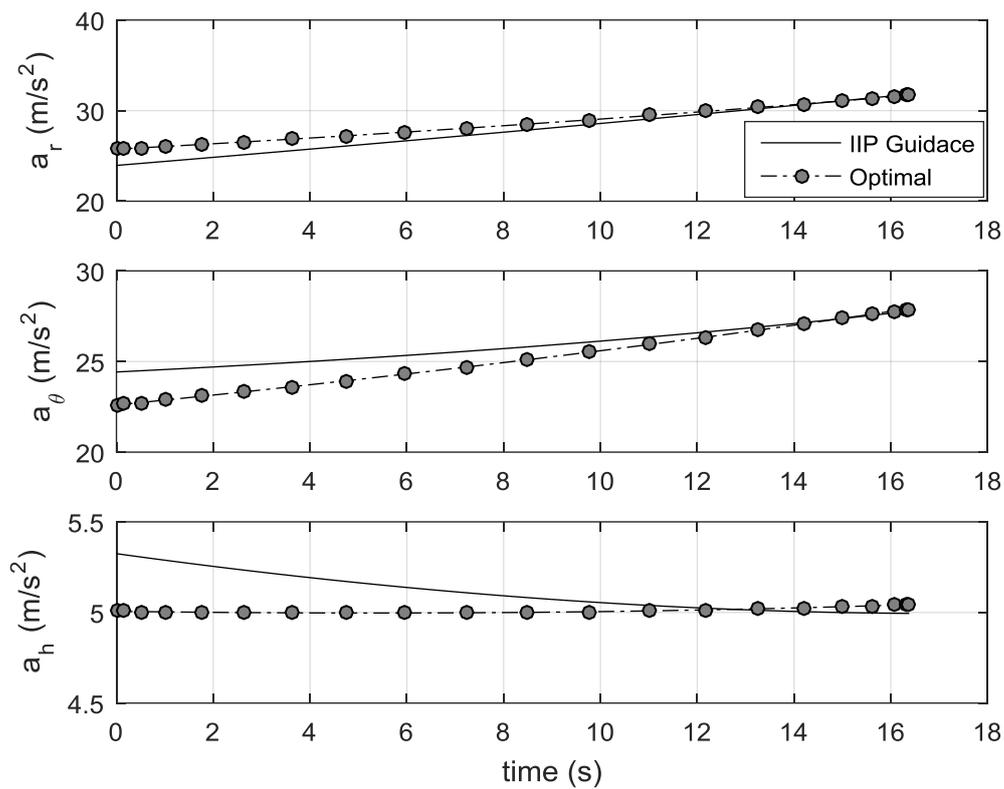

**Figure 10. Histories of acceleration components obtained by two methods (Case 4)**

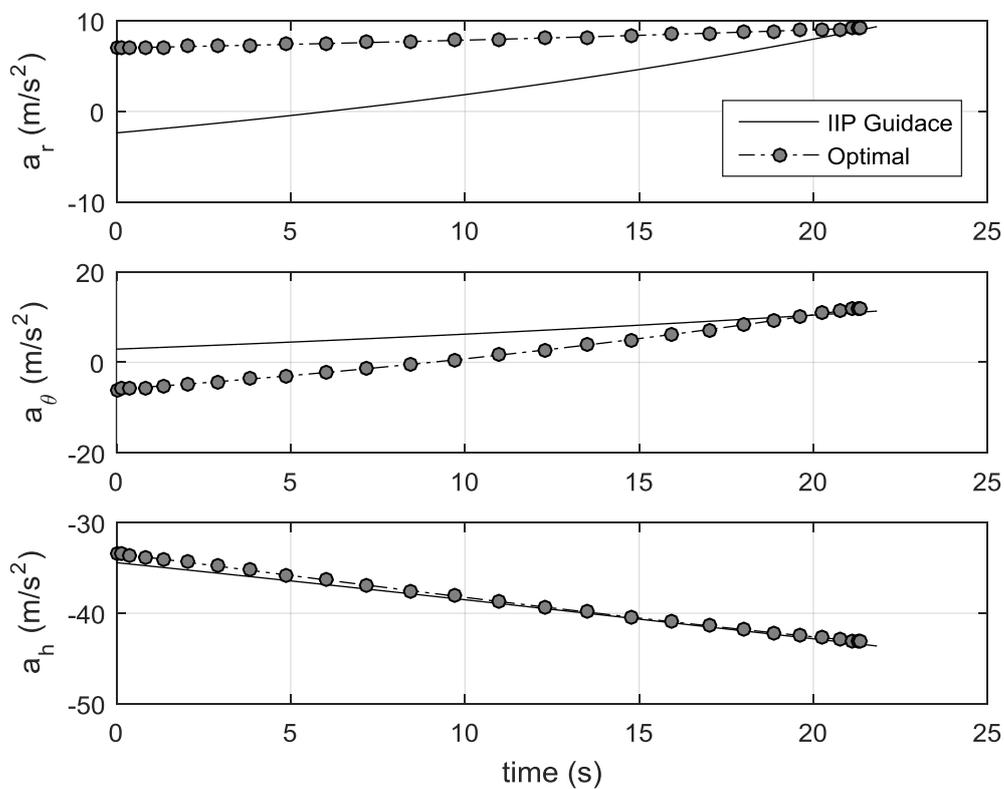

**Figure 11: Histories of acceleration components obtained by two methods (Case 5)**

While the results of the case study presented in this section demonstrate the near time-optimal performance of the proposed feedback IIP guidance law (approximately 5 % difference in final time for IIP change of 200 km), some additional analyses on the algorithm are required to justify its practicability. Additional test cases with various initial conditions and IIP change tasks should be conducted to understand the performance characteristics and robustness of the proposed algorithm. In-depth analysis on the computational load of the procedure is a potential subject for future work. Studies on the determination of time to start the proposed guidance law and the implementation of the algorithm without cut-off capability would be interesting subjects for potential future research to improve its applicability.

## V. Conclusion

A near optimal feedback guidance law to change the IIP of a rocket to a desired position was proposed in this paper. The proposed law is developed by incorporating the analytic expressions for the IIP and its time derivatives and solving a constrained optimization problem, which provides the acceleration command that aligns the IIP derivative vector with the desired direction and maximizes its magnitude. Case studies on simulation of the IIP guidance law for a rocket (separated stage of a launch vehicle) with various target points and their comparisons with numerically obtained open-loop trajectory optimization results were conducted. The results of the case study demonstrated the near time-optimal performance of the proposed guidance law. Furthermore, it can be potentially used as an on-board algorithm for moving the IIP of a rocket, after additional verification/validation and studies for performance improvement.

## Acknowledgements

This work was prepared at the Korea Advanced Institute of Science and Technology, Department of Aerospace Engineering, under a research grant from the National Research Foundation of Korea (NRF-2013M1A3A3A02042461). Authors thank the National Research Foundation of Korea for the support of this work.